# Triangle Angles and Sides in Progression and the diophantine equation $x^2+3y^2=z^2$

## 1. *Introduction*

Virtually almost every student in the late high school or early college / university years is familiar with triangles whose sidelengths are of the form $\alpha = \frac{\beta}{\sqrt{3}}$, $\beta$, and $\gamma = \frac{2\beta}{\sqrt{3}}$, where $\beta$ is a positive real number. All these triangles form a class of similar triangles with angles $A = 30°$, $B = 60°$, and $\Gamma = 90°$ . Obviously, the angles $A, B, \Gamma,$ in that order, form an arithmetic progression with difference $\varphi = 30°$. Another very familiar class of similar triangles is the set of all triangles whose sides are given by $\alpha = 3\delta$, $\beta = 4\delta$, $\gamma = 5\delta$, where $\delta$ is a positive real number. These are all right triangles, but there is another feature of theirs that stands out: their side lengths $\alpha, \beta, \gamma$ (in that order) form an arithmetic progression of difference $d = \delta$.

<u>Notation:</u> 1) *Throughout this paper, capital letters such as $A, B, \Gamma$; can either indicate the three internal angles in a well-defined triangle; or, of course, they can indicate points, such as the vertices of a triangle. For example, when we refer to the angle $A$ in a triangle $AB\Gamma$, we will always mean or refer to the internal angle in that triangle, whose vertex is $A$. In other words, the angle* $B\hat{A}\Gamma = \Gamma\hat{A}B.$

2) *Line segments will be denoted by the usual notation, e.g.* $\overline{AB}$ *denotes the straight line segment connecting or joining*

*the two points $A$ and $B$.*

3) *Sidelengths will be denoted by lower case letters. For example in a triangle $AB\Gamma$, the three sidelengths will be denoted by $\alpha, \beta, \gamma$:*

$$\alpha = \left|\overline{B\Gamma}\right| = \left|\overline{\Gamma B}\right|, \beta = \left|\overline{A\Gamma}\right| = \left|\overline{\Gamma A}\right|, \gamma = \left|\overline{AB}\right| = \left|\overline{BA}\right|$$

*Throughout this paper we will assume the standard triangle facts; namely that if $A, B, \Gamma$ are the angles of a triangle, measured in degrees, and $\alpha, \beta, \gamma$, are the side lengths, then (without loss of generality)*

$$0° < A \leq B \leq \Gamma < 180°, A + B + \Gamma = 180°,$$

*and (correspondingly) $0° < \alpha \leq \beta \leq \gamma$; and the triangle inequalities that must hold; $\alpha < \beta + \gamma, \beta < \alpha + \gamma$, and $\gamma < \alpha + \beta$ (note that if three real numbers satisfy the three triangle inequalities, they must all be positive).*

Recall that a sequence of three real numbers, $\alpha_1, \alpha_2, \alpha_3$, is an **arithmetic progression** if, and only if, $2\alpha_2 = \alpha_1 + \alpha_3$; it is a **geometric progression** if, and only if, $\alpha_2^2 = \alpha_1\alpha_3$; it is a **harmonic progression** if, and only if, the sequence $\frac{1}{\alpha_1}, \frac{1}{\alpha_2}, \frac{1}{\alpha_3}$, is an arithmetic progression.

*Throughout this paper, we refer to the corresponding sequence (of the angles of a triangle) $A, B, \Gamma$, as* ***(a)****; and the sequence (of the sidelengths) $\alpha, \beta, \gamma$, as* ***(s)****.*

There are some immediate observations:

*In order that the sequence* **(a)** *be an arithmetic progression as it follows from the* **relations** $2B=A+\Gamma$ *and* $A+B+\Gamma=180°$ *a necessary and sufficient condition is* $B=60°$*. Moreover it follows from the Law of cosines,* $\beta^2=\alpha^2+\gamma^2-2\alpha\gamma\cos B$*; that is,* $\beta^2=\alpha^2+\gamma^2-\alpha\gamma$*. Also, note that no Pythagorean Triangle (a right triangle with integer sidelengths) can have a sequence* **(a)** *which is an arithmetic progression, for the mere reason that in that case it would be* $A=30°, B=60°, \Gamma=90°$*, and thus the sidelengths* $\alpha, \beta, \gamma,$ *would satisfy,* $\alpha=\frac{\beta}{\sqrt{3}}$*,* $\beta$*, and* $\gamma=\frac{2\beta}{\sqrt{3}}$*; however, at most two of these three real numbers can be positive integers, and that happens precisely for* $\beta=m\sqrt{3}$*, where* $m$ *is a positive integer; whereas when* $\beta$ *is a positive integer, both* $\alpha$ *and* $\gamma$ *are positive irrational numbers.*

Now let us say a few words about the structure of this paper. In the next section, we present a few interesting results pertaining to the sequences **(s)** and **(a)**; in section 3 we present an analysis, of triangles whose sequences **(a)** are arithmetic progressions. In section 4 we describe the set of all triangles with integer side lengths and whose sequences **(a)** are arithmetic progressions; as we will see, to achieve that goal, we will utilize a certain subset of (the set of) parametric solutions in positive integers to the diophantine equation $x^2+3y^2=z^2$.

We finish the paper with section 5 (computations) in which we offer numerical examples (actually we offer a list of twelve such integer-sided such triangles).

To proceed further, we will need four formulas: one for the radius $r$ of a triangle's inner circle (that is, the circle inscribed in the triangle) in terms of the side lengths $\alpha, \beta, \gamma$; and the three formulas that express the inner radius $r$ in terms of the side lengths and the tangents of the half-angles $\frac{A}{2}$, $\frac{B}{2}$, and $\frac{\Gamma}{2}$.

*We denote the perimeter of a triangle by $2\tau; 2\tau = \alpha + \beta + \gamma$; so the half-perimeter is given by $\tau = \dfrac{\alpha + \beta + \gamma}{2}$.*

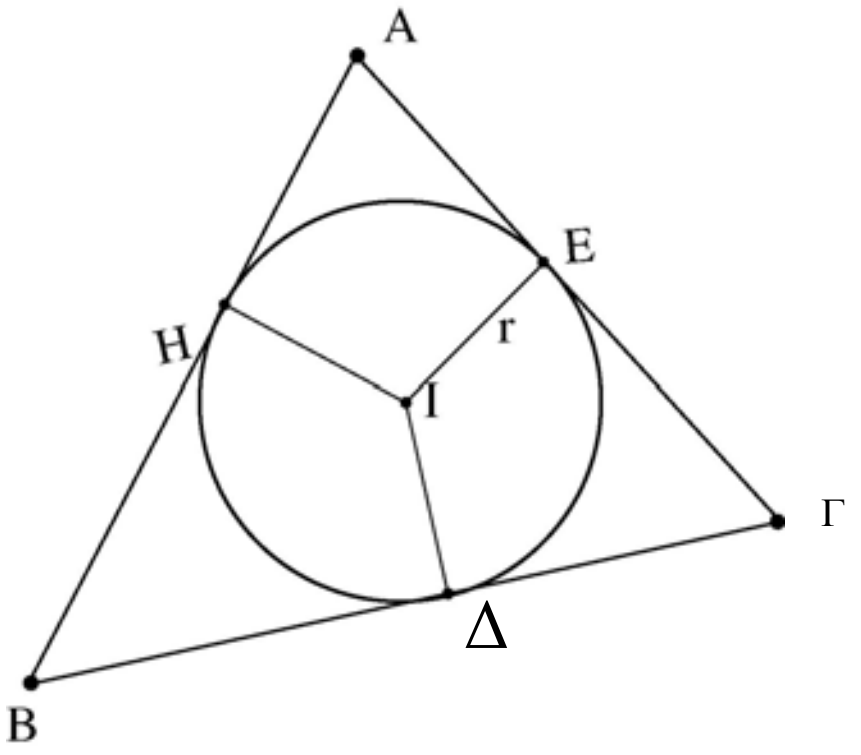

$$|\overline{AH}| = |\overline{AE}| = x \qquad r = |\overline{IE}|$$

$$|\overline{BH}| = |\overline{B\Delta}| = y \qquad r = |\overline{IH}|$$

$$|\overline{\Gamma\Delta}| = |\overline{\Gamma E}| = z \qquad r = |\overline{I\Delta}|$$

**Figure 1**

It is apparent from the picture that $2x + 2z = 2\tau; x + y + z = \tau$

$x + y = \gamma, y + z = \alpha, x + z = \beta$. And $r = y \bullet \tan\left(\dfrac{B}{2}\right) = z \bullet \tan\left(\dfrac{\Gamma}{2}\right) = x \bullet \tan\left(\dfrac{A}{2}\right)$;

$$r = (\tau - \alpha)\tan\left(\frac{A}{2}\right) = (\tau - \beta)\tan\left(\frac{B}{2}\right) = (\tau - \gamma)\tan\left(\frac{\Gamma}{2}\right) \qquad (1)$$

We employ the formula for the tangent of the sum of two angles:

$\tan(\theta_1+\theta_2)=\dfrac{\tan\theta_1+\tan\theta_2}{1-\tan\theta_1\tan\theta_2}$, which is valid for any two angles $\theta_1$ and $\theta_2$, such that neither $\theta_1$ and $\theta_2$, nor their sum $\theta_1+\theta_2$ is of the form $180°K+90°$, K an integer; we have,

$$\tan\left(\frac{A}{2}+\frac{B}{2}\right)=\frac{\tan\left(A/2\right)+\tan\left(B/2\right)}{1-\tan\left(A/2\right)\tan\left(B/2\right)} \qquad (2)$$

Also, $\tan=\left(\dfrac{A}{2}+\dfrac{B}{2}\right)=\tan\left(90°-\dfrac{\Gamma}{2}\right)=\cot\left(\dfrac{\Gamma}{2}\right)=\dfrac{1}{\tan\left(\Gamma/2\right)}$;

which in conjuction with (1) and (2) gives,

$$\frac{\dfrac{r}{\tau-\alpha}+\dfrac{r}{\tau-\beta}}{1-\dfrac{r^2}{(\tau-\alpha)(\tau-\beta)}}=\frac{\tau-\gamma}{r}\Leftrightarrow \text{ (when simplified)}$$

$$\Leftrightarrow\left[r(\tau-\alpha)+r(\tau-\beta)\right]r=(\tau-\gamma)\left[(\tau-\alpha)(\tau-\beta)-r^2\right]\Leftrightarrow$$

$$\Leftrightarrow r^2=\frac{(\tau-\alpha)(\tau-\beta)(\tau-\gamma)}{\tau} \text{ or, equivalently, } r=\sqrt{\frac{(\tau-\alpha)(\tau-\beta)(\tau-\gamma)}{\tau}} \qquad (3)$$

A Historical Perspective

*There are a few published works related to the material found in this paper. We mention these works here, starting with the most recent one. First, there is the 1997 paper of Beauregard and Suryanarayan (see reference [3]), in which Heron triangles (triangles with rational side lengths and area) with integer sidelengths and whose sidelengths are in arithmetic progression, were treated. In that paper, the diophantine equation* $x^2-3y^2=z^2$ *is used to find all such Heron triangles.*

*Then there is the 1958 short paper by A. Maxwell (see [4]), which briefly examines triangles whose angles are in arithmetic progression.*

*Further back in time, in the late nineteenth and early part of the 20th century, there are two papers by K. Schwering and J. Heinrichs. In K. Schwering's paper in 1886, he discusses triangles with integral sides and one of whose angles is double another.*

*In J. Heinrich's 1911 paper, the problem is generalized, considering the condition* $A = nB + \Gamma$, $n$ *a positive integer; between the angles* $A, B, \Gamma$ *in a triangle. Also, in a seven-page paper published in 1913, K. Schwering tackled a more general problem, by taking any linear relation between the angles. For the K. Schwering and J. Heinrich papers, see reference [5].*

## 2. *A few results*

Some of the results listed below were found by this author in an old trigonometry book, long out of print, published by a local, albeit obscure press, in Athens, Greece (see reference [2]); in that book, published in Greek only, some of the results offered below, are presented as exercises. Presumably, some of these results, one would think, could similarly be found in a smattering of other such obscure books of old, scattered around the globe. We have,

***(i)*** *In a triangle, the sequence* ***(s)*** *of the sidelengths is an arithmetic progression if, and only if, the sequence* $cot\left(A/2\right)$, $cot\left(B/2\right)$, $cot\left(\Gamma/2\right)$, *is also an arithmetic progression.*

***(ii)*** *The sequence* ***(s)*** *is an arithmetic progression if, and only if,* $tan\left(\frac{A}{2}\right)\ tan\left(\frac{\Gamma}{2}\right)= \frac{1}{3}$.

***(iii)*** *The sequence (of the squares of the side lengths)* $\alpha^2, \beta^2, \gamma^2$, *is an arithmetic progression if, and only if, the sequence* $cotA, cotB, cot\Gamma,$ *is also an arithmetic progression.*

***(iv)*** *The only triangles for which both sequences* ***(a)*** *and* ***(s)*** *are arithmetic progressions, are the equilateral triangles (the case of the trivial progressions; i.e. the difference of those progressions being zero).*

***(v)*** *The only triangles with the sequences* ***(a)*** *being arithmetic progressions and the sequences* ***(s)*** *being geometric progressions, are the equilateral ones (i.e. arithmetic progressions with difference zero and geometric progressions with ratio equal to 1).*

***(vi)*** *The only right triangles that have sequences* ***(s)*** *which are arithmetic progressions, are those whose side lengths are given by,* $\alpha=3\delta, \beta=4\delta, \gamma=5\delta; \delta$ *is a positive real number.*

***(vii)*** *The sequence of the side lenghts* $\alpha, \beta, \gamma$ *is a harmonic progression if, and only if, the sequence* $sin^2\left(\frac{A}{2}\right)$, $sin^2\left(\frac{B}{2}\right)$, $sin^2\left(\frac{\Gamma}{2}\right)$, *is also a harmonic progression.*

Below, we offer a short proof for each of parts (i) through (iii).

**Proof of (i):** We must show the equivalence,

$$2\beta = \alpha + \gamma \Leftrightarrow 2\cot\left(\frac{B}{2}\right) = \cot\left(\frac{\mathrm{A}}{2}\right) + \cot\left(\frac{\Gamma}{2}\right)$$

Using (1) we have,

$$2\cot\left(\frac{B}{2}\right) = \cot\left(\frac{\mathrm{A}}{2}\right) + \cot\left(\frac{\Gamma}{2}\right) \Leftrightarrow$$

$$\Leftrightarrow \frac{2}{\tan\left(B/2\right)} = \frac{1}{\tan\left(\mathrm{A}/2\right)} + \frac{1}{\tan\left(\Gamma/2\right)} \Leftrightarrow$$

$$\Leftrightarrow \frac{\tau - \alpha}{r} + \frac{\tau - \gamma}{r} = \frac{2\,(\tau - \beta)}{r} \Leftrightarrow \alpha + \gamma = 2\beta$$

**Proof of (ii):** We demonstrate the equivalence,

$$2\beta = \alpha + \gamma \Leftrightarrow \tan\left(\frac{\mathrm{A}}{2}\right)\tan\left(\frac{\Gamma}{2}\right) = \frac{1}{3}$$

Again we make use of (1); but also of (3):

$$\tan\left(\frac{A}{2}\right)\tan\left(\frac{\Gamma}{2}\right) = \frac{1}{3} \Leftrightarrow \left(\frac{r}{\tau - \alpha}\right)\left(\frac{r}{\tau - \gamma}\right) = \frac{1}{3} \quad \Leftrightarrow$$

$$\Leftrightarrow \frac{r^2}{(\tau-\alpha)(\tau-\gamma)} = \frac{1}{3} \Leftrightarrow \left[\frac{(\tau - \alpha)\,(\tau - \beta)\,\left(\tau - \gamma)\right)}{\tau}\right] \bullet \frac{1}{(\tau - \alpha)\,\left(\tau - \gamma\right)} = \frac{1}{3} \quad \Leftrightarrow$$

$$\Leftrightarrow \frac{\tau - \beta}{\tau} = \frac{1}{3} \quad \Leftrightarrow \quad \frac{2\tau - 2\beta}{2\tau} = \frac{1}{3} \quad \Leftrightarrow \quad \frac{(\alpha + \beta + \gamma) - 2\beta}{\alpha + \gamma + \beta} = \frac{1}{3} \quad \Leftrightarrow$$

$$\Leftrightarrow \; 3\left(\alpha + \gamma - \beta\right) = \alpha + \gamma + \beta \Leftrightarrow \alpha + \gamma = 2\beta$$

**Proof of (iii):** We must prove the equivalence,

$$2\cot \mathrm{B} = \cot \mathrm{A} + \cot \Gamma \Leftrightarrow 2\beta^2 = \alpha^2 + \gamma^2$$

We start by employing the half-angle identity $\cot\theta = \dfrac{\cot^2\left(\theta/2\right) - 1}{2\cot\left(\theta/2\right)}$,

which is valid for any angle $\theta$ which is not of the form $180° \bullet \mathrm{K}$; where $\mathrm{K}$ is an integer. We have,

$$2\cot \mathrm{B} = \cot \mathrm{A} + \cot \Gamma \qquad \Leftrightarrow$$

$$\Leftrightarrow \quad \frac{\cot^2\left(B/2\right) - 1}{\cot\left(B/2\right)} = \frac{\cot^2\left(A/2\right) - 1}{2\cot\left(A/2\right)} + \frac{\cot^2\left(\Gamma/2\right) - 1}{2\cot\left(\Gamma/2\right)} \qquad \Leftrightarrow$$

(by (1)) $$\Leftrightarrow \quad \frac{\left(\frac{(\tau - \beta)}{r}\right)^2 - 1}{\left(\frac{(\tau - \beta)}{r}\right)} = \frac{\left(\frac{(\tau - \alpha)}{r}\right)^2 - 1}{2\left(\frac{(\tau - \alpha)}{r}\right)} + \frac{\left(\frac{(\tau - \gamma)}{r}\right)^2 - 1}{2\left(\frac{(\tau - \gamma)}{r}\right)} \qquad \Leftrightarrow$$

(by using $\beta = 2\tau - (\alpha + \gamma)$ and multiplying across by $2\tau$)

$$\Leftrightarrow \quad 2\tau(\tau - \beta) - 2(\tau - \alpha)(\tau - \gamma) = \tau\beta - (\tau - \beta)\beta \qquad \Leftrightarrow$$

$$\Leftrightarrow \quad 2\tau^2 - 2\tau\beta - 2\tau^2 + 2(\alpha + \gamma)\tau - 2\alpha\gamma = \tau\beta - \tau\beta + \beta^2 \qquad \Leftrightarrow$$

$$\Leftrightarrow \quad -2\beta\tau + 2(\alpha + \gamma)\tau - 2\alpha\gamma - \beta^2 = 0 \qquad \Leftrightarrow$$

$$\Leftrightarrow \quad -2\beta\tau + 2(\alpha + \gamma)\tau - 2\alpha\gamma - \beta^2 = 0 \qquad \Leftrightarrow$$

$$\Leftrightarrow \quad -2\beta\tau + 2\tau[2\tau - \beta] - 2\alpha\gamma - \beta^2 = 0 \qquad \Leftrightarrow$$

$$\Leftrightarrow \quad -2\beta\tau + 4\tau^2 - 2\tau\beta - 2\alpha\gamma - \beta^2 = 0 \qquad \Leftrightarrow$$

$$\Leftrightarrow \quad -4\beta\tau + 4\tau^2 - 2\alpha\gamma - \beta^2 = 0 \qquad \Leftrightarrow$$

$$\Leftrightarrow \quad -2\beta(2\tau) + (2\tau)^2 - 2\alpha\gamma - \beta^2 = 0 \qquad \Leftrightarrow$$

$$\Leftrightarrow \quad -2\beta(\alpha + \beta + \gamma) + (\alpha + \beta + \gamma)^2 - 2\alpha\gamma - \beta^2 = 0 \qquad \Leftrightarrow$$

$$\Leftrightarrow \quad -2\beta\alpha - 2\beta^2 - 2\beta\gamma + \alpha^2 + \beta^2 + \gamma^2 + 2\alpha\beta + 2\beta\gamma + 2\alpha\gamma - 2\alpha\gamma - \beta^2 = 0 \Leftrightarrow$$

$$\Leftrightarrow \quad 2\beta^2 = \alpha^2 + \gamma^2$$

**Proof of (iv):** Obviously one direction in the proof is trivial: If a triangle is equilateral then $\mathrm{A} = \mathrm{B} = \Gamma = 60°$ and $\alpha = \beta = \gamma$, and so trivially, both sequences **(a)** and **(s)** are arithmetic progressions with difference equal to zero. Let us establish the converse. Namely, that if the sequences **(a)** and **(s)** are both arithmetic progressions, then the triangle must be equilateral. Indeed, by part (ii) it follows that since sequence **(s)** is an arithmetic progression, we must have, $\tan\left(\frac{A}{2}\right)\tan\left(\frac{\Gamma}{2}\right)=\frac{1}{3}$ (4)

On the other hand, since **(a)** is also an arithmetic progression, we have, in particular, $\mathrm{B}=60°$; and so $\mathrm{A}+\Gamma=120° \Rightarrow \Gamma=120°-\mathrm{A}$. Thus

(4) $\Rightarrow \tan\left(\frac{A}{2}\right)\tan\left(60° - \frac{\mathrm{A}}{2}\right)=\frac{1}{3}\Rightarrow$

$$\Rightarrow \tan\left(\frac{\mathrm{A}}{2}\right)\bullet\left[\frac{\tan 60° - \tan\left({}^{\mathrm{A}}\!/_{2}\right)}{1+\tan 60° \tan\left({}^{\mathrm{A}}\!/_{2}\right)}\right]=\frac{1}{3};$$

and since $\tan 60° = \sqrt{3}$, we obtain

$$3\tan^2\left({}^{\mathrm{A}}\!/_{2}\right) - 2\sqrt{3}\tan\left({}^{\mathrm{A}}\!/_{2}\right)+1=0 \quad \Leftrightarrow \left(\sqrt{3}\tan\left({}^{\mathrm{A}}\!/_{2}\right) - 1\right)^2=0 \;\Leftrightarrow$$

$\Leftrightarrow$ $\tan\left({}^{\mathrm{A}}\!/_{2}\right)=\frac{1}{\sqrt{3}}$ $\Leftrightarrow$ (since $0°<\frac{\mathrm{A}}{2}<90°$) $\frac{\mathrm{A}}{2}=30°;\ \mathrm{A}=60°$,

and thus, by virtue of $\mathrm{B}=60°, \Gamma=60°$.

Hence, $\mathrm{A}=\mathrm{B}=\Gamma=60°$; the triangle is equilateral.

**Proof of (v):** As with proof of part (iv), one direction is trivial. Thus, it suffices to show that if in a triangle sequence **(a)** is an arithmetic progression, while sequence **(s)** is a geometric progression, then the triangle must be equilateral. Indeed, if that were the case, we would have the simultaneous conditions $B = 60° \text{ and } \beta^2 = \alpha\gamma$. By the Law of Sines, we know that $\frac{\alpha}{\sin A} = \frac{\beta}{\sin 60°} = \frac{\gamma}{\sin \Gamma}$. Therefore,

$\beta^2 = \alpha\gamma \Rightarrow \beta^2 = \left(\frac{\beta \sin A}{\sin 60°}\right) \bullet \left(\frac{\beta \sin \Gamma}{\sin 60°}\right) \Rightarrow$ (on account of $\sin 60° = \frac{\sqrt{3}}{2}$)

$\sin A \bullet \sin \Gamma = \frac{3}{4}$. On the other hand,

$$\cos(A+\Gamma) = \cos(180° - B) = \cos 120° \Rightarrow$$

$$\Rightarrow \cos A \cos \Gamma - \sin A \sin \Gamma = -\frac{1}{2}; \quad \cos A \cos \Gamma - \frac{3}{4} = -\frac{1}{2}; \quad \cos A \cos \Gamma = \frac{1}{4}$$

We have, $\begin{Bmatrix} \sin A \sin \Gamma = \frac{3}{4} \\ \cos A \cos \Gamma = \frac{1}{4} \end{Bmatrix} \Rightarrow \begin{Bmatrix} \sin^2 A \sin^2 \Gamma = \frac{9}{16} \\ \cos^2 A \cos^2 \Gamma = \frac{1}{16} \end{Bmatrix} \Rightarrow$

$$\begin{Bmatrix} (1 - \cos^2 A)\,(1 - \cos^2 \Gamma) = \frac{9}{16} \\ \cos^2 A \cos^2 \Gamma = \frac{1}{16} \end{Bmatrix} \Rightarrow \begin{Bmatrix} \cos^2 A + \cos^2 \Gamma = \frac{1}{2} \\ \cos^2 A \cos^2 \Gamma = \frac{1}{16} \end{Bmatrix}$$

Recall that if two numbers have sum S and product P, they must be the two roots of the quadratic equation $x^2 - Sx + P = 0$. According to what we have above, the real numbers $\cos^2 A \text{ and } \cos^2 \Gamma$ must be the two roots of the quadratic equation, $x^2 - \frac{1}{2}x + \frac{1}{16} = 0 \Leftrightarrow \left(x - \frac{1}{4}\right)^2 = 0$; this quadratic equation has a double root, namely $\frac{1}{4}$. Hence,

$$\cos^2 A = \frac{1}{4} \cos^2 \Gamma \Rightarrow (\text{since } \cos A > 0, \cos \Gamma > 0)$$

$\cos A = \frac{1}{2} = \cos\Gamma \Rightarrow A = 60^\circ = \Gamma$; the triangle is equilateral.

**Proof of (vi):** In such a right triangle the sidelenghts $\alpha, \beta, \gamma$ must satisfy the two conditions,

$$\left.\begin{aligned} \gamma^2 = \alpha^2 + \beta^2 \\ 2\beta = \alpha + \gamma \end{aligned}\right\} \Rightarrow 4\gamma^2 = 4\alpha^2 + (\alpha+\gamma)^2 \Rightarrow$$

$\Rightarrow 3\gamma^2 - 2\alpha\gamma - 5\alpha^2 = 0$;

$$(3\gamma - 5\alpha)(\gamma + \alpha) = 0 \Rightarrow (\text{since } \alpha + \gamma > 0)$$

$\Rightarrow 3\gamma - 5\alpha = 0; \gamma = \frac{5\alpha}{3}$ and from $2\beta = \alpha + \gamma$ we obtain $\beta = \frac{4\alpha}{3}$;

and by setting $\alpha = 3\delta$, where $\delta$ is a positive real number we arrive at, $\alpha = 3\delta, \beta = 4\delta, \gamma = 5\delta$.

> Note that result (*vi*) has an alternate proof given in the Beauregard and Suryanarayan paper listed in [3]

**Proof of (vii):** We must prove the equivalence,

$$\frac{2}{\beta} = \frac{1}{\alpha} + \frac{1}{\gamma} \Leftrightarrow \frac{2}{\sin^2(B/2)} = \frac{1}{\sin^2(A/2)} + \frac{1}{\sin^2(\Gamma/2)}$$

We have, $\frac{2}{\sin^2(B/2)} = \frac{1}{\sin^2(A/2)} + \frac{1}{\sin^2(\Gamma/2)} \Leftrightarrow$

$$2\csc^2\left(\frac{B}{2}\right) = \csc^2\left(\frac{A}{2}\right) + \csc^2\left(\frac{\Gamma}{2}\right) \Leftrightarrow$$

$$\Leftrightarrow 2\left[1 + \cot^2(B/2)\right] = 1 + \cot^2(A/2) + 1 + \cot^2(\Gamma/2) \Leftrightarrow$$

$$\Leftrightarrow 2\cot^2(B/2) = \cot^2(A/2) + \cot^2(\Gamma/2) \Leftrightarrow$$

$$\Leftrightarrow \quad \frac{2}{\tan^2\left(B/2\right)} = \frac{1}{\tan^2\left(A/2\right)} + \frac{1}{\tan^2\left(\Gamma/2\right)} \quad \Leftrightarrow$$

$$(\text{by } (1)) \Leftrightarrow \quad \frac{2\,(\tau\text{-}\beta)^2}{r^2} = \frac{(\tau\text{-}\alpha)^2}{r^2} + \frac{(\tau\text{-}\gamma)^2}{r^2} \quad \Leftrightarrow$$

$$\Leftrightarrow 2\tau^2 - 2\tau(2\beta) + 2\beta^2 = \tau^2 - 2\tau\alpha + \alpha^2 + \tau^2 - 2\tau\gamma + \gamma^2 \Leftrightarrow$$

$$\Leftrightarrow 2\tau(\alpha + \gamma - 2\beta) = \alpha^2 + \gamma^2 - 2\beta^2 \Leftrightarrow$$

$$\Leftrightarrow (\alpha + \gamma + \beta)(\alpha + \gamma - 2\beta) = \alpha^2 + \gamma^2 - 2\beta^2 \Leftrightarrow$$

$$\Leftrightarrow \alpha^2 + \gamma^2 - 2\beta^2 + 2\alpha\gamma - \beta\gamma - \alpha\beta = \alpha^2 + \gamma^2 - 2\beta^2 \Leftrightarrow$$

$$\Leftrightarrow 2\alpha\gamma = \beta\gamma + \alpha\beta \Leftrightarrow \frac{2\,\alpha\gamma}{\alpha\beta\gamma} = \frac{\beta\gamma}{\alpha\beta\gamma} + \frac{\alpha\beta}{\alpha\beta\gamma} \Leftrightarrow \frac{2}{\beta} = \frac{1}{\alpha} + \frac{1}{\gamma}$$ ❑

## 3. *Triangle angles in progression*

It is immediately apparent that the sequence of the angles **(a)** is an arithmetic progression if, and only if, $B = 60°$ ; this follows from the conditions $2B = A + \Gamma$ and $A + B + \Gamma = 180°$.

Furthermore from the Law of Cosines we have,

$$\beta^2 = \alpha^2 + \gamma^2 - 2\alpha\gamma\cos 60° \Leftrightarrow \beta^2 = \alpha^2 + \gamma^2 - \alpha\gamma \Leftrightarrow (\alpha+\gamma)^2 = \beta^2 + 3\alpha\gamma\ .$$

In the next few lines that follow, we will see that given two positive real numbers $\beta$ and $\tau$, a unique triangle with $B = 60°$ is determined such that $0 < \alpha \le \beta \le \gamma$; as long as the ratio $\rho = \frac{2\tau}{\beta}$ falls in the semi-closed interval $(2,3]$ ; in fact the conditions $\beta > 0, \tau > 0, \text{and } 2 < \rho \le 3$ are the precise (necessary and sufficient) conditions that uniquely determine a triangle with $B = 60°$. In other words, if we are given two positive reals $\beta$ and $\tau$, such that $2 < \rho \leq 3$; where $\rho = \frac{2\tau}{\beta}$, a unique triangle can be constructed in the usual euclidean geometric way that satisfies the above conditions. Let us see why. If we denote by $S$ the sum of $\alpha$ and $\gamma$; and by $P$ their product, the relation we have above, $(\alpha+\gamma)^2 = \beta^2 + 3\alpha\gamma$ becomes $S^2 = \beta^2 + 3P \Leftrightarrow P = \frac{S^2 - \beta^2}{3}$.

On the other hand, the reals $\alpha$ and $\gamma$ are the two roots of the quadratic equation $x^2 - Sx + P = 0$. Since $\alpha$ and $\gamma$ are real roots, the discriminant of this quadratic equation must be nonnegative:

$(-S)^2 - 4P \geq 0 \Leftrightarrow S^2 - 4P \geq 0 \Leftrightarrow S^2 - 4\left(\frac{S^2 - \beta^2}{3}\right) \geq 0 \Leftrightarrow 4\beta^2 \geq S^2 \Leftrightarrow$ (since both $\beta$ and $S$ are positive) $2\beta \geq S$ and since $S = 2\tau - \beta$, we obtain $2\beta \geq 2(\tau - \beta) \Leftrightarrow 3\beta \geq 2\tau \Leftrightarrow$ (by virtue of $\beta > 0$) $\frac{3\beta}{\beta} \geq \frac{2\tau}{\beta}$; but also (obviously) $2\tau > \beta$; $\frac{2\tau}{\beta} > 1$. Thus, a necessary condition that emerges is $1 < \frac{2\tau}{\beta} \leq 3 \Leftrightarrow 1 < \rho \leq 3$.

Moreover, going back to the quadratic equation $x^2 - Sx + P = 0$, we see that since $\alpha$ and $\gamma$ are its roots with $\alpha \leq \gamma$ we must have (from the quadratic formula), $\alpha = \frac{S - \sqrt{S^2 - 4P}}{2}, \gamma = \frac{S + \sqrt{S^2 - 4P}}{2}$

Using the already familiar relations $P = \frac{S^2 - \beta^2}{3}$ and $S = 2\tau - \beta$ we obtain,

$$\alpha = \frac{2\tau - \beta - \sqrt{\frac{(3\beta - 2\tau)(\beta + 2\tau)}{3}}}{2}, \quad \gamma = \frac{2\tau - \beta + \sqrt{\frac{(3\beta - 2\tau)(\beta + 2\tau)}{3}}}{2}$$

or, equivalently,

$$\alpha = \frac{\beta\left(\frac{2\tau}{\beta} - 1\right) - \sqrt{\frac{\beta^2\left(3 - \frac{2\tau}{\beta}\right)\left(1 + \frac{2\tau}{\beta}\right)}{3}}}{2}; \quad \text{and}$$

$$\gamma = \frac{\beta\left(\frac{2\tau}{\beta} - 1\right) + \sqrt{\frac{\beta^2\left(3 - \frac{2\tau}{\beta}\right)\left(1 + \frac{2\tau}{\beta}\right)}{3}}}{2};$$

and by also using $\rho = \frac{2\tau}{\beta}$ we finally obtain,

$$\alpha = \frac{\beta}{2}\left[\rho - 1 - \sqrt{\frac{(3 - \rho)(1 + \rho)}{3}}\right], \gamma = \frac{\beta}{2}\left[\rho - 1 + \sqrt{\frac{(3 - \rho)(1 + \rho)}{3}}\right] \qquad (5)$$

But there is an additional condition: namely $\alpha>0$ which by (5) is equivalent (since $\beta > 0$) to

$$\rho - 1 - \sqrt{\frac{(3 - \rho)(1 + \rho)}{3}} > 0 \Leftrightarrow$$

$$\Leftrightarrow \rho - 1 > \sqrt{\frac{(3 - \rho)(1 + \rho)}{3}} \Leftrightarrow \text{ (under } 1 < \rho \leq 3\text{) } 3(\rho-1)^2 > (3 - \rho)(1 + \rho) \Leftrightarrow$$

$$\Leftrightarrow 3\rho^2 - 6\rho + 3 > 3 + 2\rho - \rho^2 \Leftrightarrow 4\rho^2 - 8\rho > 0 \Leftrightarrow 4\rho(\rho - 2) > 0;$$

and since we also know that $1 < \rho \leq 3$, we conclude that $\boxed{2 < \rho \leq 3}$, which is the precise condition sought after.

Recall the Law of Sines: $\frac{\alpha}{\sin A} = \frac{\beta}{\sin 60°} = \frac{\gamma}{\sin\Gamma} = 2R$ (where $R$ is the radius of the circumscribed circle); in view of $\sin B = 60° = \frac{\sqrt{3}}{2}$ we see that, $\frac{\alpha}{\sin A} = \frac{2\beta}{\sqrt{3}} = \frac{\gamma}{\sin\Gamma}$; and therefore combined with (5) we deduce that

$$\sin A = \frac{\sqrt{3}}{4}\left[\rho - 1 - \sqrt{\frac{(3 - \rho)(1 + \rho)}{3}}\right] \text{ and}$$

$$\sin \Gamma = \frac{\sqrt{3}}{4}\left[\rho - 1 + \sqrt{\frac{(3 - \rho)(1 + \rho)}{3}}\right].$$

Also,

$0 < \alpha \le \beta \le \gamma \Leftrightarrow 0° < A \le B \le \Gamma < 180° \Leftrightarrow$ (since $A + B + \Gamma = 180°$)

$0 < \sin A \le \sin B \le \sin \Gamma \le 1 \Leftrightarrow 0 < \sin A \le \sin 60° \le \sin \Gamma \le 1 \Leftrightarrow$

$0 < \sin A \le \frac{\sqrt{3}}{2} \le \sin \Gamma \le 1 \Leftrightarrow$

$$0 < \frac{\sqrt{3}}{4}\left[\rho - 1 - \sqrt{\frac{(3 - \rho)(1 + \rho)}{3}}\right] \le \frac{\sqrt{3}}{2} \le \frac{\sqrt{3}}{4}\left[\rho - 1 + \sqrt{\frac{(3 - \rho)(1 + \rho)}{3}}\right] \le 1$$

Note that if the first equal sign holds true, then so must the second one and conversely; this happens precisely when $\rho = 3$ which corresponds to the case of an equilateral triangle, i.e. $A = B = \Gamma = 60°$. Moreover in the fourth inequality, the equal sign holds true precisely for $\rho = 1 + \sqrt{3}$ (this is less apparent and it requires a little bit of algebra); this corresponds to the case of $A = 30°, B = 60°, \Gamma = 90°$. If you do the algebra, you will see that the last inequality is equivalent to $\left[\rho - (1 + \sqrt{3})\right]^2 \ge 0$.

> Conversely, by using (5), together with $2 < \rho \le 3$ and the Law of Sines, we can establish that angle $B = 60°$.

In closing, one can easily establish the three triangle inequalities $\alpha + \beta > \gamma, \alpha + \gamma > \beta, \beta + \gamma > \alpha$: the third is seen by inspection on account of $\gamma \ge \beta \ge \alpha > 0$; the second one also follows easily from $\alpha + \gamma = 2\tau - \beta$:

indeed, $\alpha+\gamma>\beta \Leftrightarrow 2\tau-\beta>\beta \Leftrightarrow 2\tau>2\beta \Leftrightarrow$ (since $\beta>0$) $\frac{2\tau}{\beta}>2$;

$\rho>2$ which is true because we have assumed that $\rho$ satisfies the key condition; namely, $2<\rho\le 3$.

We now prove triangle inequality $\alpha+\beta>\gamma$. By (5),

$\alpha+\beta>\gamma \Leftrightarrow$ (since $\beta>0$)

$$\Leftrightarrow \frac{\sqrt{3}}{4}\bullet\left[\rho+1-\sqrt{\frac{(3-\rho)(1+\rho)}{3}}\right.+1>\frac{\sqrt{3}}{4}\bullet\left[\rho+1+\sqrt{\frac{(3-\rho)(1+\rho)}{3}}\right]$$

$\Leftrightarrow$ (multiply across by 4)

$$\sqrt{3}\bullet(\rho+1)+4-\sqrt{3}\bullet(\rho+1)>2\sqrt{(3-\rho)(1+\rho)}\Leftrightarrow$$

$$\Leftrightarrow 2>\sqrt{(3-\rho)(1+\rho)}\Leftrightarrow 4>(3-\rho)(1+\rho)\Leftrightarrow$$

$$\Leftrightarrow \rho^2-2\rho+1>0\Leftrightarrow(\rho-1)^2>0,$$ which is true since $\rho\neq 1$, on account of $2<\rho\le 3$

## *Conclusions*

*In a triangle with angles* $A, B, \Gamma$*; and sidelengths* $\alpha, \beta, \gamma$*; the sequence* ***(a)*** *of the angles is an arithmetic progression if, and only if,* $B=60°$.

*In such triangle* $2<\rho \leq 3$*, where* $\rho = \frac{2\tau}{\beta}$*, and* $2\tau$ *is the triangle's perimeter.*

*Conversely, given two positive real numbers* $2\tau$ *and* $\beta$ *and such that* $2<\rho \leq 3$*, a unique triangle can be constructed that has angles* $B=60°$ *and with the two sides containing the angle* $B$ *having lengths given by*

$$\alpha=\frac{\beta}{2}\left[\rho - 1 - \sqrt{\frac{(3-\rho)(1+\rho)}{3}}\right] \text{ and } \gamma = \frac{\beta}{2}\left[\rho - 1 + \sqrt{\frac{(3-\rho)(1+\rho)}{3}}\right];$$

*and also,* $sin A = \frac{\sqrt{3}}{4}\left[\rho - 1 - \sqrt{\frac{(3-\rho)(1+\rho)}{3}}\right]$ *and*

$$sin\Gamma=\frac{\sqrt{3}}{4}\left[\rho - 1 + \sqrt{\frac{(3-\rho)(1+\rho)}{3}}\right].$$

*When* $\rho = \sqrt{3}+1$*, and only then, the triangle is a right one with* $A=30°, B=60°, \Gamma=90°$*; while precisely for* $\rho=3$*, the triangle is equilateral.*

## 4. *A determination and description of the set of all integer-sided triangles whose angles are in arithmetic progression*

As it will become evident, to be able to parametrically describe all triangles whose sidelengths are integers and whose sequences of angles **(a)** are arithmetic progressions, we will need the general solution of the three-variable diophantine equation $x^2+3y^2=z^2$. We remark here that the three-variable diophantine equations of the form $x^2+ny^2=z^2$ (where $n$ is a given positive integer), are well understood and their general solutions have been known for a long time. The interested reader should refer to [1] for further details. Without elaborating further we state the parametric formulas that describe the solution set, in positive integers, to the diophantine equation $x^2+3y^2=z^2$:

$$x=\frac{d\bullet\left|3\kappa^2-\lambda^2\right|}{2},\qquad y=d\kappa\lambda,\qquad z=\frac{d\,(3\kappa^2+\lambda^2)}{2} \qquad (6)$$

where $d,\kappa,\lambda$ are positive integers with $(\kappa,\lambda)=1$; in other words $\kappa$ and $\lambda$ are relatively prime.

Let us now consider a triangle with the sidelengths $\alpha,\beta,\gamma$ being integers and with the sequence **(a)** being an arithmetic progression; which means that $B=60°$ and consequently we obtain the by now familiar equation (from the Law of Cosines), $\beta^2=\alpha^2+\gamma^2-\alpha\gamma \Leftrightarrow \gamma^2-\alpha\gamma+\alpha^2-\beta^2=0$ . Keeping in mind that $0<\alpha\le\beta\le\gamma$ and solving the last quadratic equation for $\gamma$ we obtain,

$$\gamma=\frac{\alpha\pm\sqrt{4\beta^2-3\alpha^2}}{2} \qquad (7)$$

In what follows, we make use of a well-known result in elementary number theory. According to the result, if n and m are positive integers; then $\sqrt[n]{m}$ (the nth root of $m$) is a rational number if, and only if $m$ is an nth integer power; that is $m = k^n$, for some positive integer $k$. Which means that $\sqrt[n]{m} = k$, a positive integer. In particular, the square root of a positive integer is rational if, and only if, that positive integer is an integral (or perfect) square. The interested reader who maybe unfamiliar with the above described result, may wish to consult either of the two references; [6] or [7]. The above comes into play because according to (7); (since $\beta, \gamma, \alpha$ are integers) the square root $\sqrt{4\beta^2 - \alpha^2}$ must be a rational number.

We easily see that γ will be a positive integer if, and only if, $4\beta^2 - 3\alpha^2 = \delta^2$, where $\delta$ is a positive integer; $\frac{\alpha \pm \delta}{2} > 0$, and the integers $\alpha$ and $\delta$ have the same parity (they are both odd or both even). Obviously $\frac{\alpha + \delta}{2} > 0$, since both $\alpha$ and $\delta$ are positive integers.

Note however, that the possibility $\gamma = \frac{\alpha - \delta}{2}$ is contradictory. This is easy to see: for, on account of $\gamma > 0$, it would require $\alpha > \delta$; however if that were the case, then from $4\beta^2 - 3\alpha^2 = \delta^2 \Rightarrow 4\beta^2 < 3\alpha^2 + \alpha^2 \Rightarrow 4\beta^2 < 4\alpha^2 \Rightarrow$ (since both $\beta$ and $\alpha$ are positive) $\beta < \alpha$, contrary to $\alpha \le \beta$. Hence we conclude that in (7), the plus sign must hold:

$$\left.\begin{aligned} \gamma = \frac{\alpha + \sqrt{4\beta^2 - 3\alpha^2}}{2} \\ 4\beta^2 = 3\alpha^2 + \delta^2 \\ \text{with } \alpha, \beta, \gamma, \delta \in \mathbb{Z}^+ \end{aligned}\right\} \Leftrightarrow \left.\begin{aligned} \gamma = \frac{\alpha + \delta}{2} \\ (2\beta)^2 = \delta^2 + 3\alpha^2 \\ \text{with } \alpha, \beta, \gamma, \delta \in \mathbb{Z}^+ \end{aligned}\right\} \qquad (8)$$

(set of the positive integers)

The second equation in (8) shows that the triple $(\delta, \alpha, 2\beta)$ is a positive integer solution to the diophantine equation $x^2 + 3y^2 = z^2$: we must have $\delta = \frac{d \bullet |3\kappa^2 - \lambda^2|}{2}$, $\alpha = d\kappa\lambda$, $2\beta = \frac{d\,(3\kappa^2 + \lambda^2)}{2}$, which when combined with (8) produces,

$$\left\{\begin{aligned} \alpha = d\kappa\lambda,\ \beta = \frac{d\,(3\kappa^2 + \lambda^2)}{4},\ \gamma = d\left(\frac{2\kappa\lambda + |3\kappa^2 - \lambda^2|}{4}\right) \\ \text{where } d, \kappa, \lambda \in \mathbb{Z}^+, \text{ and } (\kappa, \lambda) = 1 \end{aligned}\right\} \qquad (9)$$

First observe that the integer parameters $\kappa, \lambda$ must either be both odd or one of them even, the other odd; this, because of the condition $(\kappa, \lambda) = 1$. If $\kappa$ and $\lambda$ have different parity, then $d$ must be a multiple of 4, in order to obtain integer values for $\beta$ and $\gamma$; if on the other hand, $\kappa \equiv \lambda \equiv 1$ (mod 2); then $\kappa^2 \equiv \lambda^2 \equiv 1$ (mod 4) and so, $3\kappa^2 + \lambda^2 \equiv 0$ (mod 4), which shows that $d$ can take any positive integer value; $\beta$ and $\gamma$ will be positive integers regardless. Also note that the three integers $i\kappa\lambda$, $i\left(\frac{3\kappa^2 + \lambda^2}{4}\right)$, $i\left(\frac{2\kappa\lambda + |3\kappa^2 - \lambda^2|}{4}\right)$, where $i = 1$ (when $\kappa \equiv \lambda \equiv 1$ (mod 2)) or $i = 4$ (when $\kappa + \lambda \equiv 1$ (mod 2)) are mutually relatively prime (i.e. any two of them will be relatively prime), as long as $\lambda$ is not divisible by 3; if $\lambda$ is a multiple of 3, then the greatest common divisor of any two of those three integers will equal 3; we leave this as an exercise for the interested reader to prove. Next, we impose the condition $\alpha \leq \beta \leq \delta$ on formulas (9). This will, in turn, enable us to determine the conditions which the integer parameters $\kappa \text{ and } \lambda$ must satisfy. We have,

Insofar, we have obtained all the positive integer solutions to the equation $\beta^2 = \gamma^2 + \alpha^2 - \alpha\gamma$; next we must determine those solutions among them that satisfy $\alpha \leq \beta \leq \gamma$. We apply (9):

$$d\kappa\lambda \leq \frac{d\,(3\kappa^2 + \lambda^2)}{4} \leq d\left(\frac{2\kappa\lambda + \left|3\kappa^2 - \lambda^2\right|}{4}\right) \qquad (10)$$

Since $d, \kappa, \lambda$ are all positive, a bit of algebra shows that the first inequality in (10) is equivalent to,

$$\left(\frac{\lambda}{\kappa} - 1\right)\left(\frac{\lambda}{\kappa} - 3\right) \geq 0 \Leftrightarrow \left(\frac{\lambda}{\kappa} \leq 1 \text{ or } \frac{\lambda}{\kappa} \geq 3\right) \qquad (11)$$

Note that the cases $\frac{\lambda}{\kappa} = 1$ and $\frac{\lambda}{\kappa} = 3$ correspond to $(\lambda = \kappa = 1)$ and $(\lambda = 3, \kappa = 1)$ , in view of $(\lambda, \kappa) = 1$, and by (9) we see that in the first case, $\alpha = \beta = \gamma = d$; while in the second, $\alpha = \beta = \gamma = 3d$; thus in both cases the triangle obtained (assuming there is one, as it is shown below) would be equilateral.

Since λ and κ are positive integers, (11) can also be expressed as $(1 \leq \lambda \leq \kappa \text{ or } 3 \leq 3\kappa \leq \lambda)$; it can be expressed in interval notation as well: $\frac{\lambda}{\kappa} \varepsilon (0,1] \cup [3 + \infty)$. Next we show that under (11), the second inequality in (10) is satisfied as well. Note that $3\kappa^2 - \lambda^2$ cannot be zero since $\sqrt{3}$ is an irrational number. If $3\kappa^2-\lambda^2>0$, then $\left|3\kappa^2-\lambda^2\right|=3\kappa^2-\lambda^2$ and so the second inequality in (10) becomes $d\left(\frac{2\kappa\lambda+3\kappa^2-\lambda^2}{4}\right)\geq\frac{d(3\kappa^2+\lambda^2)}{4}\Leftrightarrow$ $\left(\frac{\lambda}{\kappa}\right)\left[\frac{\lambda}{\kappa}-1\right]\leq 0\Leftrightarrow$ (since $0<\frac{\lambda}{\kappa}$), $0<\frac{\lambda}{\kappa}\leq 1$, which is consistent with (11) and also with $3\kappa^2-\lambda^2>0$ (i.e. $0<\frac{\lambda}{\kappa}<\sqrt{3}$). Likewise if $3\kappa^2-\lambda^2<0$, then $\left|3\kappa^2-\lambda^2\right|=\lambda^2-3\kappa^2$ and the second inequality is equivalent to $\frac{\lambda}{\kappa}\geq 3$ which is consistent with (11) and also with $3\kappa^2-\lambda^2<0$ (i.e. $\sqrt{3}<\frac{\lambda}{\kappa}$). It is now clear that (9) together with (11) describe all positive integer solutions to $\beta^2=\alpha^2+\gamma^2-2\alpha\gamma$ such that $\alpha\leq\beta\leq\gamma$. Moreover, the Law of Cosines ensures us that triangle angle $B$ is equal to $60°$.

The only thing left to check is the triangle inequalities, $\alpha+\beta>\gamma$, $\beta+\gamma>\alpha$, $\alpha+\gamma>\beta$, to make sure that an actual triangle is formed:

We have $\alpha+\beta>\gamma\Leftrightarrow$ (by (9)) $2\kappa\lambda+3\kappa^2+\lambda^2>\left|3\kappa^2-\lambda^2\right|$, which is obviously true since $\left|3\kappa^2-\lambda^2\right|\leq\left|3\kappa^2\right|+\left|\lambda^2\right|=3\kappa^2+\lambda^2$, and $\kappa,\lambda$ are positive. The second triangle inequality, $\beta+\gamma>\alpha$, is equivalent to $(\kappa-\lambda)^2+2\kappa^2+\left|3\kappa^2-\lambda^2\right|>0$; obviously true. To establish the third triangle inequality $\alpha+\gamma>\beta$, one needs to make use of (11) and to distinguish between the cases $3\kappa^2-\lambda^2>0$ and $3\kappa^2-\lambda^2<0$. We omit the details. Also, observe that if we combine (9) with (11) and do some algebra, we will see that the ratio $\rho=\frac{2\tau}{\beta}=\frac{\alpha+\beta+\gamma}{\gamma}$, indeed satisfies $2<\rho\leq 3$, as expected from the conclusion of the previous section.

## *Conclusions*

*The parametric description of all triangles with integer sidelengths $\alpha, \beta, \gamma$ and such that $\alpha \le \beta \le \gamma$, and whose sequence **(a)** of angles is an arithmetic progression is as follows:*

$$\alpha = d\kappa\lambda,\ \beta = \frac{d\left(3\kappa^2+\lambda^2\right)}{4},\ \gamma = d\left(\frac{2\kappa\lambda+\left|3\kappa^2-\lambda^2\right|}{4}\right),$$

*where $d, \kappa, \lambda$ are positive integer parameters with $(\kappa,\lambda)=1$ and with $(1 \le \lambda \le \kappa \text{ or } 3\kappa \le \lambda)$ (or equivalently, in interval notation, $\frac{\lambda}{\kappa} \varepsilon (0,1] \cup [3+\infty)$. When $\kappa, \lambda$ are both odd, then $d$ can be any positive integer; while for $\kappa+\lambda \equiv 1$ (mod 2), $d$ must be a multiple of 4. Furthermore, angle $B=60°$ and (from the Law of Sines)*

$$\sin A \frac{\alpha \sin 60°}{\beta} = \frac{2\,\kappa\lambda\sqrt{3}}{3\kappa^2+\lambda^2},$$

*with (necessarily) $0<A \le 60°$; and with $\varphi = 60° - A$, where $\varphi$ is the difference of the arithmetic progression **(a);** $0° \le \varphi < 60°$.*

*Finally, angle $\Gamma = 120° - A = 2\varphi + A = \varphi + 60°$. Also, note that for a given pair $(\kappa,\lambda)$ by varying $d$ an entire class of similar triangles is generated.*

*Moreover, as expected, the ratio $\rho = \frac{2\tau}{\beta}$ falls in the interval $(2,3]$; $2 < \rho \le 3$.*

## 5. *Computations*

Below we list all triangles with integer sidelengths and whose sequences **(a)** are arithmetic progressions, and with the parameters $\lambda$ and $\kappa$ satisfying (in addition to the other two conditions) $1 \le \lambda, \kappa \le 5$; and the parameter $d$ being equal to 4 when $\lambda$, $\kappa$ having different parity; while $d = 1$, when $\lambda$, $\kappa$ are both odd. The data you will find below, are calculated from the formulas of the previous section, and with aid of a scientific calculator for approximately determining the value of the angle $A$. There are twelve triangles listed below.

> Note: The examples below show that, although an entire class of similar triangles is obtained by fixing $\lambda$ and $\kappa$ and letting $d$ vary, these classes are not disjoint. This is due to the fact that $\lambda$ is allowed to be a multiple of 3, in which case a triangle is obtained, with each of its sides being a multiple of 3.

**1)** $\kappa = \lambda = 1, d = 1, \alpha = \beta = \gamma = 1, \rho = 3, \sin A = \frac{\sqrt{3}}{2},$
$A = 60°, \varphi = 0°, B = 60°, \Gamma = 60°$

**2)** $\kappa = 2, \lambda = 1, d = 4, \alpha = 8, \beta = 13, \gamma = 15, \rho = \frac{36}{13} \approx 2.769230769$,

$\sin A = \frac{4\sqrt{3}}{13}$, $A \approx 32.2042275°, \varphi \approx 27.7957725°, B=60°, \Gamma \approx 87.7957725°$

**3)** $\kappa = 3, \lambda = 1, d = 1, \alpha = 3, \beta = 7, \gamma = 8, \rho = \frac{72}{28} = \frac{18}{7} \approx 2.57142851,$

$\sin A = \frac{3\sqrt{3}}{14}, A \approx 21.7867893°, \varphi \approx 38.2132107°, B = 60°, \Gamma \approx 98.2132107°$

**4)** $\kappa = 4, \lambda = 1, d = 4, \alpha = 16, \beta = 49, \gamma = 55, \rho = \frac{120}{49} \approx 2.448979592,$

$\sin A = \frac{8\sqrt{3}}{49}, A \approx 16.4264214°, \varphi \approx 43.5735786°, B = 60°, \Gamma \approx 103.5735786°$

**5)** $\kappa = 5, \lambda = 1, d = 1, \alpha = 5, \beta = 19, \gamma = 21,$

$\rho = \frac{45}{19} \approx 2.368421053, \sin A = \frac{10\sqrt{3}}{76} = \frac{5\sqrt{3}}{38},$

$A \approx 13.17355111°, \varphi \approx 46.82644889°, B = 60°, \Gamma \approx 106.8264489°$

**6)** $\kappa = 3, \lambda = 2, d = 4, \alpha = 24, \beta = 31, \gamma = 35,$

$\rho = \frac{90}{31} \approx 2.903225806, \sin A = \frac{12\sqrt{3}}{31},$

$A \approx 42.10344887°, \varphi \approx 17.89655113°, B = 60°, \Gamma \approx 77.89655113°$

**7)** $\kappa = 5, \lambda = 2, d = 4, \alpha = 40, \beta = 79, \gamma = 91, \rho = \frac{210}{79} \approx 2.658227848,$

$\sin A = \frac{20\sqrt{3}}{79}, A \approx 26.00782389°, \varphi \approx 33.99217611°, B = 60°, \Gamma \approx 93.99217611°$

**8)** $\kappa = 1, \lambda = 3, d = 1, \alpha = \beta = \gamma = 3, \sin A = \frac{\sqrt{3}}{2}, A = 60°, \varphi = 0°,$

$B = 60°, \Gamma = 60°$

**9)** $\kappa = 4, \lambda = 3, d = 4, \alpha = 48, \beta = 57, \gamma = 63,$

$\rho = \frac{168}{57} = \frac{56}{19} \approx 2.947368421, \sin A = \frac{24\sqrt{3}}{57} = \frac{8\sqrt{3}}{19},$

$A \approx 46.82644889°, \varphi \approx 13.17355111°, B = 60°, \Gamma \approx 73.17355111°$

**10)** $\kappa = 5, \lambda = 3, d = 1, \alpha = 15, \beta = 21, \gamma = 24,$

$$\rho = \frac{60}{21} = \frac{20}{7} \approx 2.857142857, \ \sin A = \frac{30\sqrt{3}}{84} = \frac{5\sqrt{3}}{14},$$
$$A \approx 38.2132107°, \ \varphi \approx 21.7867893°, \ B = 60°, \ \Gamma \approx 81.7867893°$$

**11)** $\kappa = 1, \lambda = 4, d = 4, \alpha = 16, \beta = 19, \gamma = 21,$

$$\rho = \frac{56}{19} \approx 2.94736841, \ \sin A = \frac{32\sqrt{3}}{76} = \frac{8\sqrt{3}}{19},$$
$$A \approx 46.82644889°, \ \varphi \approx 13.17355111°, \ B = 60°, \ \Gamma \approx 73.17355111°$$

Note the similarity of this triangle with the triangle in #9 $(\kappa = 4, \lambda = 3, d = 4)$

**12)** $\kappa = 1, \lambda = 5, d = 1, \alpha = 5, \beta = 7, \gamma = 8,$

$$\rho = \frac{80}{28} = \frac{20}{7} \approx 2.857142857, \ \sin A = \frac{5\sqrt{3}}{14},$$
$$A \approx 38.2132107°, \ \varphi \approx 21.7867893°, \ B = 60°, \ \Gamma \approx 81.7867893°$$

Note the similarity of this triangle with the triangle in #10 $(\kappa = 5, \lambda = 3, d = 1)$